\documentclass[11pt]{article}  
\begin{document} \newtheorem{prop}{Proposition}[section]
\newtheorem{Def}{Definition}[section] \newtheorem{theorem}{Theorem}[section]
\newtheorem{lemma}{Lemma}[section] \newtheorem{Cor}{Corollary}[section]

\title{\bf Global solutions with infinite energy for the 1-dimensional Zakharov system}
\author{{\bf
 Hartmut Pecher}\\
 Fachbereich Mathematik und Naturwissenschaften\\
  Bergische Universit\"at Wuppertal\\
 Gau{\ss}str.  20
  \\ D-42097 Wuppertal\\
   Germany\\
    e-mail Hartmut.Pecher@math.uni-wuppertal.de}
     \date{}
 \maketitle

\begin{abstract}
The 1-dimensional Zakharov system is shown to have a unique global solution for data without finite 
energy. The proof uses the " I-method " introduced by Colliander, Keel, Staffilani, Takaoka, and Tao in 
connection with a refined bilinear Strichartz estimate.
\end{abstract}

\setcounter{section}{-1}
\section{Introduction} Consider the (1+1)-dimensional Cauchy problem for the Zakharov system
\begin{eqnarray}
\label{1}
 iu_t + u_{xx} & = & nu
\\ \label{2}
n_{tt}-n_{xx} & = & (|u|^2)_{xx}
\\
\label{3}
 u(0) \quad = \quad u_0 \quad , \quad
n(0) \,& = &\, n_0 \quad , \quad n_t(0) \quad = \quad n_1
 \end{eqnarray}
where $u$ is a complex-valued und $n$ a
real-valued function defined for $(x,t) \in {\bf R}\times{\bf R}^+$.

The Zakharov system was introduced in \cite{Z} to describe Langmuir turbulence in a plasma.

Our main result is the existence of a unique global solution for data without finite energy, more 
precisely
we assume $ u_0\in H^s({\bf R})$ , $ n_0 \in H^{s-1}({\bf R}) $ , $ A^{-1/2}n_1 \in H^{s-1}({\bf R}) $ ,
where $ 1>s>5/6 $ , $ A:=-\frac{d^2}{dx^2}$.

This result can be proven by using the conservation laws, namely conservation of $\|u(t)\| $ and $$
E(u,n):=\|u_x(t)\|^2 + \frac{1}{2}(\|n(t)\|^2 + \|A^{-1/2}n_t(t)\|^2) + \int_{-\infty}^{\infty}
n(t)|u(t)|^2\, dx $$ although under our assumptions these quantities are not finite, in general.

Results of this type were given in various situations in the last years in the framework of the Fourier
restriction norm method in most of the applications. One approach is to use Bourgain's trick to split
the data into high and low frequency parts. He used it to prove global well-posedness for the (2+1)- 
and (3+1)-dimensional Schr\"odinger equations with rough data without finite energy \cite{B1},\cite{B2}
and for the wave equation \cite{B3}. Later it was also used for other model equations 
\cite{CST},\cite{FLP},\cite{KT},\cite{KPV},\cite{P2}. Concerning the problem at hand the author had been 
able to show global 
well-posedness for data $(u_0,n_0,n_1)\in H^s \times L^2 \times \dot{H}^{-1} $ for $1>s>9/10$ \cite{P1}. 
Remark here that no data $n_0 \not\in L^2$ were admissible because in such a case the nonlinear part of 
$n(t)$ could not be shown to belong to $L^2$ which is necessary for this method. In contrast, the approach 
here allows data $n_0 \not\in L^2$ and also less regular data $u_0$.

Another approach was initiated by Colliander, Keel, Staffilani, Takaoka and Tao in \cite{CKSTT3}, called 
the I-method. The main idea is to use a modified energy functional which is also defined for less regular 
functions and not strictly conserved. When one is able to control its growth in time explicitly this 
allows 
to iterate a modified local existence theorem to continue the solution to any time $T$ and moreover to 
estimate its growth in time. This method was successfully applied by these authors to several equations 
which have a scaling invariance with sometimes even optimal global well-posedness results. It was 
used in \cite{CKSTT3} to improve Bourgain's global well-posedness results \cite{B1},\cite{B2} for the 
(2+1)- and (3+1)-dimensional Schr\"odinger equation with a further improvement in \cite{CKSTT6}. Later it 
was applied to the (1+1)-dimensional derivative Schr\"odinger equation \cite{CKSTT1} with an (almost) 
optimal result in \cite{CKSTT2} and to the KdV and modified KdV equation with also optimal results in some 
cases \cite{CKSTT4},\cite{CKSTT5}.

Although in our situation such a scaling argument does not work we are able to suitably modify the method 
to prove the above mentioned global existence result for the Zakharov system.

The paper is organized as follows. We transform the system in the usual way into a first order system. 
Then 
we apply the multiplier $I_N$ for given $s<1$ and $N>>1$ to it, where $ \widehat{I_N 
f}(\xi):=m_N(\xi)\widehat{f}(\xi) $. Here $m_N(\xi)$ is a smooth, radially symmetric, and nonincreasing 
function of $|\xi|$, defined by $ m_N(\xi) = 1 $ for $| \xi|\le N$ and $m_N(\xi) = (\frac{N}{|\xi|})^{1-s} 
$ for $ |\xi| \ge 2N $. We drop $N$ from the notation for short and remark that $I:H^s \to H^1$ is a 
smoothing operator in the following sense:
$$ \|u\|_{X^{m,b}_{\varphi}} \le c \|Iu\|_{X^{m+1-s,b}_{\varphi}} \le c N^{1-s}\|u\|_{X^{m,b}_{\varphi}} 
$$
Here we used the $X^{m,b}_{\varphi}$ - spaces which are defined as follows: for an equation of the form 
$iu_t-\varphi(-i\partial_x)u = 0$, where $\varphi$ is a measurable function, let $X^{m,b}_{\varphi}$ be 
the 
completion of ${\cal S}({\bf R}\times{\bf R})$ with respect to
\begin{eqnarray*}
\|f\|_{X^{m,b}_{\varphi}} & := & \|\langle\xi\rangle^m\langle\tau\rangle^b {\cal 
F}(e^{it\varphi(-i\partial_x)}f(x,t))\|_{L^2_{\xi 
\tau}} \\ & = & \|\langle\xi\rangle^m\langle\tau+\varphi(\xi)\rangle^b \widehat{f}(\xi,\tau)\|_{L^2_{\xi 
\tau}} 
\end{eqnarray*}
For $\varphi(\xi) = \pm | \xi|$ we use the notation $X^{m,b}_{\pm}$ and for $\varphi(\xi) = |\xi|^2$ 
simply $X^{m,b}$. For a given time interval $I$ we define $\|f\|_{X^{m,b}(I)} = \inf_{\tilde{f}_{|I} = f} 
\|\tilde{f}\|_{X^{m,b}} $ and similarly $\|f\|_{X^{m,b}_{\pm}(I)} $.

For the modified (by $I$ multiplied) Zakharov system we then prove a local existence theorem by using the 
precise estimates given by \cite{GTV} for the standard Zakharov system in connection with an interpolation 
type lemma in \cite{CKSTT5}. Our aim is to extract a factor $T^{\delta}$ with maximal $\delta$ from the 
nonlinear estimates in order to give an optimal lower bound for the local existence time $T$ in terms of 
the norms of the data. Because the difference of the differentiability classes of the data is maximal 
(=1), one is forced also to use here the auxiliary spaces $Y^m_{\varphi}$ (cf. \cite{GTV}), defined by
\begin{eqnarray*}
\|f\|_{Y^m_{\varphi}}: & = & \|\langle\xi\rangle^m \langle\tau\rangle^{-1} {\cal 
F}(e^{-it\varphi(-i\partial_x)}f(x,t))\|_{L^2_{\xi}L^1_{\tau}} \\
& = & \|\langle\xi\rangle^m\langle\tau +\varphi(\xi)\rangle^{-1} 
\widehat{f}(\xi,\tau)\|_{L^2_{\xi}L^1_{\tau}}
\end{eqnarray*}
As is typical for the $I$-method one then has to consider in detail the modified energy functional 
$E(Iu,In)$ and to control its growth in time in dependence of the time interval and the parameter $N$ (cf. 
the definition of $I$ above). The increment of the energy has to be small for small time intervals and 
large $N$. Because the modified energy functional is somehow close to the original one here some sort of 
cancellation helps. An important tool is also a refined Strichartz estimate for the product of a wave and 
a Schr\"odinger part along the lines of Bourgain's improvements for the simpler pure Schr\"odinger case 
(cf. Lemma \ref{Lemma A}). This estimate for the modified energy functional can also control the growth of 
the corresponding norms of the solution of the problem during its time evolution. One iterates the local 
existence theorem with time steps of equal length in order to reach any given fixed time $T$. To achieve 
this one has to make the process uniform which can be done if $s$ is close enough to $1$ (namely $s>5/6$).

We collect some elementary facts about the spaces $X^{m,b}_{\varphi}$ and $Y^m_{\varphi}$.\\ 
The following interpolation property is well-known:\\
$ X^{(1-\Theta)m_0+\Theta m_1,(1-\Theta)b_0+\Theta b_1}_{\varphi} = 
(X^{m_0,b_0}_{\varphi},X^{m_1,b_1}_{\varphi})_{[\Theta]} $ for $ \Theta \in [0,1] $. \\
If $u$ is a 
solution of $iu_t+\varphi(-i\partial_x)u = 0$ with $u(0)=f$ and $\psi$ is a cutoff function in 
$C^{\infty}_0({\bf R})$ with $supp \, \psi \subset (-2,2)$ , $\psi \equiv 1$ on $[-1,1]$ , $ \psi(t) = 
\psi(-t) $ , $ \psi(t)\ge 0 $ , $\psi_{\delta}(t):=\psi(\frac{t}{\delta}) \, , $  $ 0<\delta \le 1$, we 
have 
for $b>0$: 
$$\|\psi_1 u\|_{X^{m,b}_{\varphi}} \le c \|f\|_{H^m} $$
If $v$ is a solution of the problem $iv_t +\varphi(-i\partial_x)v = F $ , $ v(0)=0 $ , we have for $b'+1 
\ge b \ge 0 \ge b' > -1/2$
$$ \|\psi_{\delta} v \|_{X^{m,b}_{\varphi}} \le c \delta^{1+b'-b} \|F\|_{X^{m,b'}_{\varphi}} $$
and, if $b'+1 \ge b \ge 0 \ge b'$, we have
$$ \|\psi_{\delta} v \|_{X^{m,b}_{\varphi}} \le c (\delta^{1+b'-b} \|F\|_{X^{m,b'}_{\varphi}} + 
\delta^{\frac{1}{2}-b} \|F\|_{Y^m_{\varphi}}) $$
(for a proof cf. \cite{GTV}, Lemma 2.1). Moreover, if $w(t) = \int_0^t e^{i(t-s)\varphi(-i\partial_x)} 
F(s) \, ds $ we have by \cite{GTV}, Lemma 2.2, especially (2.35), for $\delta \le 1 $
\begin{equation}
\label{b}
\|w\|_{C^0([0,\delta],H_x^1)} \le c \|F\|_{Y^1_{\varphi}[0,\delta]}
\end{equation}
Finally, if $1/2 > b > b' \ge 0$ , $ m \in {\bf R} $ , we have the embedding
\begin{equation}
\|f\|_{X^{m,b'}_{\varphi}[0,\delta]} \le c \delta^{b-b'} \|f\|_{X^{m,b}_{\varphi}[0,\delta]}
\end{equation}
For the convenience of the reader we repeat the proof of \cite{G}, Lemma 1.10. The claimed estimate is an 
immediate consequence of the following
\begin{lemma}
For $1/2 > b > b' \ge 0$ , $ 0 < \delta \le 1 $ , $ m \in {\bf R} $ the following estimate holds:
$$ \| \psi_{\delta} f \|_{X^{m,b'}_{\varphi}} \le c \delta^{b-b'} \|f\|_{X^{m,b}_{\varphi}} $$
\end{lemma}
{\bf Proof:} The following Sobolev multiplication rule holds:
$$ \|fg\|_{H^{b'}_t} \le c \|f\|_{H^{\frac{1}{2}-(b-b')}_t} \|g\|_{H^b_t} $$
This rule follows easily by the Leibniz rule for fractional derivatives, using $J^s := {\cal 
F}^{-1}\langle\tau\rangle^s {\cal F}$:
\begin{eqnarray*}
\|fg\|_{H^{b'}_t} & \le & c(\|(J^{b'}f)g\|_{L^2_t} + \|f(J^{b'}g)\|_{L^2_t}) \\
& \le & c(\|J^{b'}f\|_{L^p_t} \|g\|_{L^{p'}_t} + \|f\|_{L^{q'}_t} \|J^{b'}g\|_{L^q_t})
\end{eqnarray*}
with $\frac{1}{p}=b$ , $ \frac{1}{p'} = \frac{1}{2}-b $ , $ \frac{1}{q'}=b-b' $ , $ \frac{1}{q} = 
\frac{1}{2} - (b-b') $. Sobolev's embedding theorem gives the claimed result. Consequently we get
$$ \|\psi_{\delta}g\|_{H^{b'}_t} \le c \|\psi_{\delta}\|_{H^{\frac{1}{2}-(b-b')}_t} \|g\|_{H^b_t} \le c 
\delta^{b-b'} \|g\|_{H^b_t} $$
and thus
\begin{eqnarray*}
\|\psi_{\delta}f\|_{X^{m,b'}_{\varphi}} = \|e^{it\varphi(-i\partial_x)}\psi_{\delta}f\|_{H^m_x\otimes 
H^{b'}_t} & \le & c \delta^{b-b'}\|e^{it\varphi(-i\partial_x)}f\|_{H^m_x \otimes H^b_t}\\ 
& = & c\delta^{b-b'}\|f\|_{X^{m,b}_{\varphi}} 
\end{eqnarray*}
Fundamental are the following linear Strichartz type estimates for the Schr\"odinger equation (cf. e.g. 
\cite{GTV}, Lemma 2.4):
$$ \|e^{it\partial^2_x} \psi\|_{L^q_t(I,L_x^r({\bf R}))} \le c \|\psi\|_{L^2_x({\bf R})} $$
and
$$ \|u\|_{L^q_t(I,L_x^r({\bf R}))} \le c \|u\|_{X^{0,\frac{1}{2}+}(I)} $$
if $ 0 \le \frac{2}{q} = \frac{1}{2} - \frac{1}{r} $ , especially
$$ \|u\|_{L^6_{xt}} \le c \|u\|_{X^{0,\frac{1}{2}+}} $$
which by interpolation with the trivial case $\|u\|_{L^2_{xt}} = \|u\|_{X^{0,0}}$ gives:
$$ \|u\|_{L^p_{xt}} \le c \|u\|_{X^{0,\frac{3}{2}(\frac{1}{2}-\frac{1}{p})+}} $$
if $ 2 < p \le 6 $. For the wave equation we only use $ \|n_{\pm}\|_{L^{\infty}_t L^2_x} \le c 
\|n_{\pm}\|_{X^{0,\frac{1}{2}+}_{\pm}}$.\\
We use the notation $\langle \lambda \rangle := (1+\lambda^2)^{1/2} $. Let $a \pm$ denote a number 
slightly larger (resp., smaller) than $a$.\\
{\bf Acknowledgement:} I thank Axel Gr\"unrock for very helpful discussions.

\section{Local existence}
The system (\ref{1}),(\ref{2}),(\ref{3}) has the following conserved quantities:
$$ \|u(t)\|  =: M $$
and
$$ E(u,n):= \|A^{1/2}u(t)\|^2 + 1/2(\|n(t)\|^2 + \|V(t)\|^2) + \int_{-\infty}^{+\infty} n(t)|u(t)|^2 \, dx 
$$
where $ V_x := -n_t$ and $ A:=-\frac{d^2}{dx^2}$ . \\
The system (\ref{1}),(\ref{2}),(\ref{3}) is now transformed into a first order system in $t$ as follows: 
with $ n_{\pm}:= n \pm iA^{-1/2}n_t $ , i.e. $ n=\frac{1}{2}(n_++n_-) ,$ $2iA^{-1/2}n_t=n_+-n_-$, and 
$\overline{n_+} = n_- $ this gives
\begin{eqnarray}
\label{1'}
iu_t + u_{xx} & = & \frac{1}{2}(n_+ + n_-)u \\
\label{2'}
in_{\pm t} \mp A^{1/2}n_{\pm} & = & \pm A^{1/2}(|u|^2) \\
\label{3'}
u(0) = u_0 \quad , \quad n_{\pm}(0) & = & n_{\pm 0} \quad := \quad n_0 \pm iA^{-1/2}n_1
\end{eqnarray} 
The energy is given by
$$ E(u,n_+)  = \|A^{1/2}u\|^2 + \frac{1}{2} \|n_+\|^2 + \frac{1}{2} \int (n_++\overline{n_+})|u|^2 \, dx 
$$
By Gagliardo-Nirenberg
\begin{eqnarray*}
\int n|u|^2 \, dx \le \frac{1}{4}\int n^2 \, dx + c \int |u|^4 \, dx & \le & \frac{1}{4} \|n\|^2 + c 
\|u_x\| \|u\|^3 \\
& \le & \frac{1}{4}(\|n\|^2 + \|u_x\|^2) + c_0 \|u\|^6
\end{eqnarray*}
This easily implies
\begin{equation}
\label{*}
\|A^{1/2}u\|^2 + \|n\|^2 + \|V\|^2 \le c(E+\|u\|^6) = c_0(E+M^6)
\end{equation}
and also
\begin{equation}
\label{*'}
E \le c_0(\|A^{1/2}u\|^2 + \|n\|^2 + \|V\|^2 +M^6)
\end{equation}
We want to apply the I-method (for the definition of I see the introduction). A crucial role is played by 
the modified energy $E(Iu,In_+)$ for the system
\begin{eqnarray}
\label{1''}
iIu_t + Iu_{xx} & = & \frac{1}{2} I[(n_+ +n_-)u] \\
\label{2''}
iIn_{\pm t} \mp A^{1/2}In_{\pm} & = & \pm IA^{1/2}(|u|^2) \\
\label{3''}
Iu(0) = Iu_0 \, , \, In_{\pm} (0) & = & In_{\pm 0} = I(n_0 \pm iA^{-1/2}n_1)
\end{eqnarray}
namely
$$ E(Iu,In_+) := \|Iu_x\|^2 + \frac{1}{2} \|In_+\|^2 + \frac{1}{2} \int I(n_+ + \overline{n_+})|Iu|^2 \, 
dx $$
 which is not conserved but its growth is controllable.

An elementary but lengthy calculation shows
\begin{eqnarray}
\nonumber
\frac{d}{dt} E(Iu,In_+) & = & Re \langle I(n_+ + \overline{n_+}) Iu - I((n_+ + 
\overline{n_+})u),Iu_t\rangle \\
& & + \, Re \langle In_+,iA^{1/2}(|Iu|^2 - I(|u|^2))\rangle
\label{4}
\end{eqnarray}
If $I= id $ this again shows the conservation of $E(u,n_+)$.

Before considering this modified energy in detail we give a local existence result for the system 
(\ref{1''}),(\ref{2''}),(\ref{3''}), which essentially uses the bilinear estimates given by \cite{GTV} for 
their local existence result of the Zakharov system.
\begin{prop}
\label{Prop. 1}
Assume $ s \ge 1/2 $. Let $(u_0,n_{+0},n_{-0}) \in H^s \times H^{s-1} \times H^{s-1} $ be given. Then 
there exists a positive number $\delta \sim \frac{1}{(\|Iu_0\|_{H^1} + \|In_{+0}\|_{L^2} + 
\|In_{-0}\|_{L^2})^{4+}} $ such that the system (\ref{1''}), (\ref{2''}), (\ref{3''}) has a unique local 
solution in the time interval $[0,\delta]$ with the property (dropping from now on $[0,\delta]$ from the 
notation):
$$ \|Iu\|_{X^{1,\frac{1}{2}}} + \|In_+\|_{X^{0,\frac{1}{2}+}_+} + \|In_-\|_{X^{0,\frac{1}{2}+}_-} \le 
c(\|Iu_0\|_{H^1} + \|In_{+0}\|_{L^2} + \|In_{-0}\|_{L^2})$$
This solution also belongs to $C^0([0,\delta],H^1_x({\bf R}))$ and
$$ \|Iu\|_{C^0([0,\delta],H^1_x({\bf R}))} \le c ( \|Iu_0\|_{H^1} +  \|In_{+0}\|_{L^2} + \|In_{-0}\|_{L^2}) 
$$
\end{prop}
{\bf Proof:} We use the corresponding integral equations to define a mapping $S=(S_0,S_1)$ by
\begin{eqnarray*}
S_0(Iu(t)) & = & Ie^{it\partial^2_x}u_0 + \frac{1}{2}\int_0^t 
e^{i(t-s)\partial^2_x}I(u(s)((n_+(s)+n_-(s))\, ds \\
S_1(In_{\pm}(t)) & = & Ie^{itA^{1/2}}n_{\pm 0} \pm i \int_0^t e^{\mp i(t-s)A^{1/2}}A^{1/2}I(|u(s)|^2)\, ds
\end{eqnarray*}
We use \cite{GTV}, Lemma 4.3 to conclude for $s \ge 1/2 $:
$$ \|n_{\pm}u\|_{X^{s,-\frac{1}{2}}} \le c \|n_{\pm}\|_{X^{s-1,\frac{3}{8}+}_{\pm}} 
\|u\|_{X^{s,\frac{3}{8}+}} $$
Similarly \cite{GTV}, Lemma 4.5 shows
$$ \|n_{\pm}u\|_{Y^s} \le c \|n_{\pm}\|_{X^{s-1,\frac{3}{8}+}_{\pm}} \|u\|_{X^{s,\frac{3}{8}+}} $$
Finally, \cite{GTV}, Lemma 4.4 shows for $s \ge 0$:
$$\|A^{1/2}(|u|^2)\|_{X^{s-1,-\frac{1}{2}+}_{\pm}} \le c \|u\|_{X^{s,\frac{1}{4}++}}^2 $$
These estimates imply similar estimates including the I-operator by the interpolation lemma of 
\cite{CKSTT5}, namely
$$ \|I(n_{\pm}u)\|_{X^{1,-\frac{1}{2}}} + \|I(n_{\pm}u)\|_{Y^1} \le c 
\|In_{\pm}\|_{X^{0,\frac{3}{8}+}_{\pm}} \|Iu\|_{X^{1,\frac{3}{8}+}} $$
and 
$$ \|IA^{1/2}(|u|^2)\|_{X^{0,-\frac{1}{2}+}_{\pm}} \le c \|Iu\|_{X^{1,\frac{1}{4}++}}^2 $$
where $c$ is independent of $N$.\\
The same estimates also hold true for functions defined on $[0,\delta]$ only, and for such functions we 
can also use the embedding (5). This gives
$$ \|I(n_{\pm}u)\|_{X^{1,-\frac{1}{2}}} + \|I(n_{\pm}u)\|_{Y^1} \le c 
\|In_{\pm}\|_{X^{0,\frac{1}{2}+}_{\pm}} \|Iu\|_{X^{1,\frac{1}{2}}} \delta^{\frac{1}{4}-} $$
and
$$ \|IA^{1/2}(|u|^2)\|_{X^{0,-\frac{1}{2}+}_{\pm}} \le c \|Iu\|_{X^{1,\frac{1}{2}}}^2 
\delta^{\frac{1}{2}-} $$
Using these estimates the integral equations lead to (remark here that one needs the space $Y^1$, cf. 
\cite{GTV}, Lemma 2.1):
\begin{eqnarray*}
 \|S_0(Iu)\|_{X^{1,\frac{1}{2}}} & \le & c \|Iu_0\|_{H^1}+ c (\|In_+\|_{X^{0,\frac{1}{2}+}_+} + 
\|In_-\|_{X^{0,\frac{1}{2}+}_-}) \|Iu\|_{X^{1,\frac{1}{2}}} \delta^{\frac{1}{4}-} \\
 \|S_1(In_{\pm})\|_{X^{0,\frac{1}{2}+}_{\pm}} & \le & c \|In_{\pm 0}\|_{L^2} + c 
\|Iu\|_{X^{1,\frac{1}{2}}}^2 \delta^{\frac{1}{2}-}
 \end{eqnarray*}
 The standard contraction argument gives the existence of a unique solution on $[0,\delta]$ with
 $$ \|Iu\|_{X^{1,\frac{1}{2}}} + \|In_+\|_{X^{0,\frac{1}{2}+}_+} + \|In_-\|_{X^{0,\frac{1}{2}+}_-}
  \le 2c(\|Iu_0\|_{H^1} + \|In_{+0}\|_{L^2} + \|In_{-0}\|_{L^2}) $$
 provided
 $$ c \delta^{\frac{1}{4}-}(\|Iu_0\|_{H^1} + \|In_{+0}\|_{L^2}  + \|In_{-0}\|_{L^2}) <1 $$
Concerning the property $Iu \in C^0([0,\delta],H^1_x)$ we refer to \cite{GTV}, Lemma 2.2 (use the first 
integral equation and $I(un_{\pm}) \in Y^1$). Moreover (\ref{b}) gives
\begin{eqnarray*}
\|Iu\|_{C^0([0,\delta],H^1)} & \le & \|Iu_0\|_{H^1} + c (\|I(n_+u)\|_{Y^1} + \|I(n_-u)\|_{Y^1}) \\
& \le & \|Iu_0\|_{H^1} + c \delta^{\frac{1}{4}-} (\|In_+\|_{X^{0,\frac{1}{2}+}_+} + 
\|In_-\|_{X^{0,\frac{1}{2}+}_-}) \|Iu\|_{X^{1,\frac{1}{2}}} \\
& \le & \|Iu_0\|_{H^1} + c 
\frac{(\|Iu_0\|_{H^1}+\|In_{+0}\|_{L^2}+\|In_{-0}\|_{L^2})^2}{\|Iu_0\|_{H^1}+\|In_{+0}\|_{L^2}+\|In_{-0}\|
_{L^2}} \\
& \le & c(\|Iu_0\|_{H^1}+\|In_{+0}\|_{L^2}+\|In_{-0}\|_{L^2})
\end{eqnarray*}

\section{A bilinear Strichartz estimate}
\begin{lemma}
\label{Lemma 2.1}
$$ \|(D^{1/2}_x u) n_{\pm}\|_{L^2_{xt}} \le c \|n_{\pm}\|_{X^{0,\frac{1}{2}+}_{\pm}} 
\|u\|_{X^{0,\frac{1}{2}+}} $$
\end{lemma}
{\bf Proof:} We split the domain of integration into the following parts \\
a) $ supp \, \widehat{u} \subset \{ |\xi|\ge 2 \} $ \\
b) $ supp \, \widehat{u} \subset \{ |\xi| \le 2 \} $ \\
a) We assume $ supp \, \widehat{m_{\pm}} \subset \{ \xi \ge 0 \} $ (the other part $ supp \, 
\widehat{m_{\pm}} \subset \{ \xi \le 0 \} $ can be treated similarly) and $ supp \, \widehat{v} \subset \{ 
|\xi| \ge 2 \} $. In this region we conclude as follows:
\begin{eqnarray} \label{***}
& & \|e^{it\partial_x^2}D^{1/2}_x v e^{\pm it|\partial_x|} m_{\pm}\|_{L^2_{xt}}^2 \\ \nonumber
& = & \int d\xi dt \left| \int_{ \xi = \xi_1 + \xi_2 \, , \, \xi_2 \ge 0 \, , \, |\xi_1| \ge 
2}e^{-it\xi_1^2\pm it|\xi_2|} \widehat{v}(\xi_1)\widehat{m_{\pm}}(\xi_2)|\xi_1|^{\frac{1}{2}} d\xi_1 
\right|^2 \\ \nonumber
& = & \int d\xi dt \int_{\xi=\xi_1+\xi_2=\eta_1+\eta_2 \, , \, \eta_2,\xi_2\ge 0 \, , \,  
|\xi_1|,|\eta_1|\ge 2 } e^{-it(\xi_1^2\pm |\xi_2| - \eta_1^2 \mp |\eta_2|)} 
\widehat{v}(\xi_1)\overline{\widehat{v}(\eta_1)}\cdot \\ \nonumber
& & \hspace{6cm}\cdot \,   \widehat{m_{\pm}}(\xi_2) \overline{\widehat{m_{\pm}}(\eta_2)} 
|\xi_1|^{\frac{1}{2}} |\eta_1|^{\frac{1}{2}}d\xi_1 d\eta_1
 \\ \nonumber
& = & \int d\xi \int d\xi_1 d\eta_1 \delta 
(P(\eta_1))\widehat{v}(\xi_1)\overline{\widehat{v}(\eta_1)}\widehat{m_{\pm}}(\xi_2)\overline{\widehat{m_{\pm}
}(\eta_2)} |\xi_1|^{\frac{1}{2}}|\eta_1|^{\frac{1}{2}}
\end{eqnarray}
with
\begin{eqnarray*}
P(\eta_1) & := & \xi_1^2 \pm |\xi_2| - \eta_1^2 \mp |\xi - \eta_1| = \xi_1^2 \pm |\xi_2| - \eta_1^2 \mp 
|\eta_2| \\
& = & \xi_1^2 \pm \xi_2 - \eta_1^2 \mp \eta_2 = \xi_1^2 \pm (\xi -\xi_1) - \eta_1^2 \mp (\xi - \eta_1) \\
& = & \xi_1^2 - \eta_1^2 \mp (\xi_1 - \eta_1) = (\xi_1 - \eta_1)[(\xi_1+\eta_1) \mp 1]
\end{eqnarray*}
This function has the simple zeroes $\eta_1 = \xi_1$ and $\eta_1 = \pm 1 - \xi_1 $. Moreover $ P'(\eta_1) 
= -2\eta_1 \pm 1$ , thus $ |P'(\eta_1)| \sim |\eta_1| $ in our region $|\eta_1| \ge 2$. Using the 
well-known identity
$$ \int \delta(P(\eta_1))f(\eta_1)\, d\eta_1 = \sum \frac{f(x_k)}{|P'(x_k)|} $$
where $x_k$ denotes the simple zeroes of $P$, we remark that in our case for the zeroes we have $|\eta_1| 
\sim |\xi_1| $, and therefore the factor $|\xi_1|^{\frac{1}{2}} |\eta_1|^{\frac{1}{2}}$ cancels with 
$|P'(x_k)|$. Thus we can estimate (\ref{***}) using Schwarz' inequality by
\begin{eqnarray*}
& & c \int d\xi \int d\xi_1 |\widehat{v}(\xi_1)\overline{{\widehat{v}(\xi_1)}}\widehat{m_{\pm}}(\xi - 
\xi_1) \overline{\widehat{m_{\pm}}(\xi - \xi_1)}| \\
& + & c\int d\xi \int d\xi_1 |\widehat{v}(\xi_1)\overline{{\widehat{v}(\pm 1 
-\xi_1)}}\widehat{m_{\pm}}(\xi - \xi_1) \overline{\widehat{m_{\pm}}(\xi -(\pm 1 - \xi_1))}| \\
& \le & c \|\widehat{v}\|_{L^2}^2 \|\widehat{m_{\pm}}\|_{L^2}^2 = c \|v\|_{L^2}^2 \|m_{\pm}\|_{L^2}^2
\end{eqnarray*}
and the claimed estimate follows directly in the region $supp \, \widehat{u} \subset \{|\xi|\ge 2\}$ (cf. 
e.g. \cite{P2}, Lemma 1.4, \cite{G}, Lemma 2.1 or \cite{KS}, Section 3).\\
b) In the region $ supp \, \widehat{u} \subset \{ |\xi| \le 2 \} $ we have
\begin{eqnarray*}
\|(D_x^{\frac{1}{2}}u)n_{\pm}\|_{L^2_{xt}} & \le & \|D_x^{\frac{1}{2}}u\|_{L^2_t L^{\infty}_x} 
\|n_{\pm}\|_{L^{\infty}_t L^2_x} \le c \|D^{\frac{1}{2}}_x u\|_{L^2_t H^{\frac{1}{2}+}_x} 
\|n_{\pm}\|_{X^{0,\frac{1}{2}+}_{\pm}} \\
& \le & c \|u\|_{L^2_t L^2_x} \|n_{\pm}\|_{X_{\pm}^{0,\frac{1}{2}+}} \le  c \|u\|_{X^{0,\frac{1}{2}+}} 
\|n_{\pm}\|_{X_{\pm}^{0,\frac{1}{2}+}}
\end{eqnarray*}

\begin{lemma}
\label{Lemma A}
$$ \|(D^{1/2}_x u) n_{\pm}\|_{L^2_{xt}} \le c \|n_{\pm}\|_{X^{0,\frac{1}{2}+}_{\pm}} 
\|u\|_{X^{0+,\frac{1}{2}}} $$
\end{lemma}
{\bf Proof:} Interpolate the estimate of the previous lemma with
\begin{eqnarray*}
\|(D^{\frac{1}{2}}_x u)n_{\pm}\|_{L^2_{xt}} & \le & \|n_{\pm}\|_{L^{\infty}_t L^2_x} \|D^{\frac{1}{2}}_x 
u\|_{L^2_t L^{\infty}_x} \le c \|n_{\pm}\|_{X^{0,\frac{1}{2}+}_{\pm}} \|u\|_{L^2_t H^{1+}_x} \\ & \le &  c 
 \|n_{\pm}\|_{X^{0,\frac{1}{2}+}_{\pm}} \|u\|_{X^{1+,0}} 
\end{eqnarray*}

\begin{lemma}
\label{Lemma A'}
$$ \|(D^{1/2}_x u) n_{\pm}\|_{L^{2+}_t L^2_x} \le c \|n_{\pm}\|_{X^{0,\frac{1}{2}}_{\pm}} 
\|u\|_{X^{0+,\frac{1}{2}}} $$
\end{lemma}
{\bf Proof:} By Sobolev's embedding theorem and Strichartz' estimate we have
\begin{eqnarray*}
\|(D_x^{\frac{1}{2}}u)n_{\pm}\|_{L^{4-}_tL^2_x} \le c \|n_{\pm}\|_{L^{\infty -}_tL^2_x} \|D^{\frac{1}{2}}_x 
u\|_{L^4_tL^{\infty}_x} & \le & c \|n_{\pm}\|_{X^{0,\frac{1}{2}-}_{\pm}} \|D^{\frac{1}{2}}_x 
u\|_{L^4_tH_x^{\frac{1}{4}+,4}} \\
& \le & c \|n_{\pm}\|_{X^{0,\frac{1}{2}-}_{\pm}} \|u\|_{X^{\frac{3}{4}+,\frac{3}{8}+}} 
\end{eqnarray*}
Interpolation with Lemma \ref{Lemma 2.1} gives the claimed result.

A variant of this lemma is given in the following
\begin{lemma}
\label{Lemma B}
$$ \|(\widehat{D^{1/2}_x u}) \ast \widehat{n_{\pm}} \|_{L^2_{\xi}L^{2-}_{\tau}} \le c 
\|n_{\pm}\|_{X^{0,\frac{1}{2}+}_{\pm}} \|u\|_{X^{0+,\frac{1}{2}}} $$
\end{lemma}
{\bf Proof:} On one hand Lemma \ref{Lemma 2.1} gives
\begin{equation}
\label{**''}
\|(\widehat{D^{1/2}_x u}) \ast \widehat{n_{\pm}} \|_{L^2_{\xi \tau}} \le c 
\|n_{\pm}\|_{X^{0,\frac{1}{2}+}_{\pm}} \|u\|_{X^{0,\frac{1}{2}+}}
\end{equation}
On the other hand Young's inequality shows
\begin{equation}
\label{*''}
\|(\widehat{D^{1/2}_x u}) \ast \widehat{n_{\pm}} \|_{L^2_{\xi}L^{\frac{4}{3}}_{\tau}} \le c 
\|\widehat{n_{\pm}}\|_{L^2_{\xi}L^1_{\tau}}
\|\widehat{D^{1/2}_x u}\|_{L^1_{\xi} L^{\frac{4}{3}}_{\tau}}
\end{equation}
Now by Schwarz' inequality
\begin{eqnarray*}
 \|\widehat{n_{\pm}}\|_{L^2_{\xi}L^1_{\tau}} & = & \|\int |\widehat{n_{\pm}}(\xi,\tau)| \langle \tau \pm 
|\xi| \rangle ^{\frac{1}{2}+} \langle \tau \pm |\xi| \rangle ^{-\frac{1}{2}-} d\tau \|_{L^2_{\xi}} \\   & 
\le & c \|\widehat{n_{\pm}}(\xi,\tau) \langle \tau \pm |\xi| \rangle ^{\frac{1}{2}+} \|_{L^2_{\xi \tau}}  
  =  c \|n_{\pm}\|_{X^{0,\frac{1}{2}+}} 
 \end{eqnarray*}
 and by H\"older's inequality in $\tau$ and Schwarz' inequality in $\xi$:
 \begin{eqnarray*}
 \|\widehat{D^{1/2}_x u}\|_{L^1_{\xi}L^{\frac{4}{3}}_{\tau}} & = & \|\widehat{D^{1/2}_x u}(\xi,\tau) 
\langle \tau + \xi^2 \rangle ^{\frac{1}{4}+} \langle \tau + \xi^2 \rangle 
^{-\frac{1}{4}-}\|_{L^1_{\xi}L^{\frac{4}{3}}_{\tau}} \\
 & \le & c \|\widehat{D^{1/2}_x u}(\xi,\tau) \langle \tau + \xi^2 \rangle ^{\frac{1}{4}+} 
\|_{L^1_{\xi}L^2_{\tau}} \\
 & = &  c \|\widehat{D^{1/2}_x u}(\xi,\tau) \langle \tau + \xi^2 \rangle ^{\frac{1}{4}+} \langle \xi 
\rangle ^{\frac{1}{2}+}  \langle \xi \rangle ^{-\frac{1}{2}-}\|_{L^1_{\xi}L^2_{\tau}} \\
& \le & c \|\widehat{D^{1/2}_x u}(\xi,\tau) \langle \tau + \xi^2 \rangle ^{\frac{1}{4}+} \langle \xi 
\rangle ^{\frac{1}{2}+}  \|_{L^2_{\xi \tau}} \\
& \le & c \|u\|_{X^{1+,\frac{1}{4}+}}
\end{eqnarray*}
Interpolating (\ref{*''}) and (\ref{**''}) we get the result.

We also need a bilinear Strichartz' refinement for the pure Schr\"odinger problem. We have the well-known
\begin{lemma}
\label{Lemma C}
If $u_1$,$u_2$ fulfill $|\xi_1| >> |\xi_2| \ge 1$ for $\xi_i \in supp \, \widehat{u_i} $ $(i=1,2)$, the 
following estimate holds:
$$\|(D_x^{\frac{1}{2}}u_1) u_2\|_{L^2_{xt}} \le c \|u_1\|_{X^{0,\frac{1}{2}+}} 
\|u_2\|_{X^{0,\frac{1}{2}+}} $$
\end{lemma}
{\bf Proof:} \cite{CKSTT1}, Lemma 7.1.

We also have the following variant:
\begin{lemma}
\label{Lemma D}
Under the assumptions of Lemma \ref{Lemma C} we have:
$$\|(D_x^{\frac{1}{2}}u_1) u_2\|_{L^{2+}_t L^2_x} \le c \|u_1\|_{X^{0+,\frac{1}{2}}} 
\|u_2\|_{X^{0,\frac{1}{2}}} $$
\end{lemma}
{\bf Proof:} similarly as the proof of Lemma \ref{Lemma A'}.\\
{\bf Remark:} All the estimates in this section remain true, if any of the functions on the l.h.s. of the 
estimates are replaced by their complex conjugates.

\section{Estimates for the modified energy}
The main step towards global existence is an exact control of the increment of the modified energy.
\begin{prop}
\label{Prop. 3.1}
Let $(u,n_{\pm})$ be a solution of (\ref{1'}),(\ref{2'}),(\ref{3'}) on $[0,\delta]$ in the sense of Prop. 
\ref{Prop. 1}. Then the following estimate holds (for $N \ge 1 $ , $ s > 3/4 $):
\begin{eqnarray*}
\lefteqn{ |E(Iu(\delta),In_+(\delta)) - E(Iu(0),In_+(0))| } \\
&  \le &  c( (N^{-\frac{1}{2}+} \delta^{\frac{1}{2}-} + N^{-\frac{3}{2}+} \delta^{0+}) 
\|In_+\|_{X^{0,\frac{1}{2}+}_+}   \|Iu\|_{X^{1,\frac{1}{2}}}^2 \\ & &
 + (N^{-3+} + N^{-1+}\delta^{\frac{1}{2}-})\|In_+\|_{X^{0,\frac{1}{2}+}_+} ^2\|Iu\|_{X^{1,\frac{1}{2}}}^2)
\end{eqnarray*}
\end{prop}
{\bf Proof:} Using (\ref{4}) and replacing $Iu_t$ by (\ref{1''}) we have to show
\begin{equation}
\label{5}
\left|\int_0^{\delta} \int In_+ A^{1/2}(|Iu|^2 - I(|u|^2)) dxdt \right| \le 
\frac{c}{N^{1-}}\delta^{\frac{1}{2}-} \|In_+\|_{X^{0,\frac{1}{2}+}_+} \|Iu\|_{X^{1,\frac{1}{2}}}^2
\end{equation}
and
\begin{eqnarray}
\label{6} \lefteqn{
\left| \int_0^{\delta} \int (Iu)_{xx}(I(n_+u)-In_+Iu)dxdt\right| }  \\ & \le & 
c(N^{-\frac{1}{2}+}\delta^{\frac{1}{2}-} + N^{-\frac{3}{2}+} \delta^{0+}) \|In_+\|_{X^{0,\frac{1}{2}+}_+} 
\|Iu\|_{X^{1,\frac{1}{2}}}^2 \nonumber
\end{eqnarray}
as well as
\begin{equation}
\label{7}
\left| \int_0^{\delta} \int I(n_+u)(I(n_+u)-In_+Iu)dxdt\right| \le c (\frac{1}{N^{3-}} + N^{-1+} 
\delta^{\frac{1}{2}-}) \|In_+\|_{X^{0,\frac{1}{2}+}_+}^2 \|Iu\|_{X^{1,\frac{1}{2}}}^2
\end{equation}
Here and in the sequel we assume w.l.o.g. the Fourier transforms of all these functions to be nonnegative, 
ignore the appearance of complex conjugates, use dyadic decompositions w. r. to the frequencies $|\xi_j| 
\sim N_j = 2^k $ $(k=0,1,2,...)$. In order to sum over the dyadic pieces at the end we need to have extra 
factors $N_j^{0-}$ everywhere.

We start with (\ref{5}) which follows from
\begin{eqnarray} \nonumber
\lefteqn{\hspace{-2cm} \int_0^{\delta} \int_* \widehat{n_+}(\xi_1,t)|\xi_2 + \xi_3|^{\frac{1}{2}} \left| 
\frac{m(\xi_2+\xi_3)-m(\xi_2)m(\xi_3)}{m(\xi_2)m(\xi_3)} \right| 
\widehat{u_2}(\xi_2,t)\widehat{u_3}(\xi_3,t) d\xi dt } \\
& \le & \frac{c}{N^{1-}} \delta^{\frac{1}{2}-} \|n_+\|_{X^{0,\frac{1}{2}+}_+}  
\|u_2\|_{X^{1,\frac{1}{2}}}\|u_3\|_{X^{1,\frac{1}{2}}} \label{8}
\end{eqnarray}
Here and in the sequel * denotes integration over the set $\sum_{i=1}^3 \xi_i = 0 $ (or $\sum_{i=1}^4 
\xi_i = 0 $).\\
The symmetry in $\xi_2,\xi_3$ allows to assume $N_2 \ge N_3$, and moreover we can assume $N_2 \ge N$, 
because otherwise the symbol is $\equiv 0$. The condition $\sum_{i=1}^3 \xi_i =0$ implies $N_1 \le cN_2$. 
Thus $N_2 \sim N_{max}$, where $ N_{max} := \max(N_1,N_2,N_3) $.\\
We have $ |\frac{m(\xi_2+\xi_3)-m(\xi_2)m(\xi_3)}{m(\xi_2)m(\xi_3)}| \le \frac{c}{|m(\xi_2)||m(\xi_3)|} \le 
c\left\langle(\frac{N_2}{N})^{1/2}\right\rangle $ and thus the bound 
\begin{eqnarray}
\nonumber
\lefteqn{ c \int_0^{\delta} \int_* \widehat{n_+}(\xi_1,t)|\xi_2|^{1/2} 
\widehat{u_2}(\xi_2,t)\widehat{u_3}(\xi_3,t)d\xi dt \left\langle(\frac{N_2}{N})^{1/2}\right\rangle } \\
\nonumber
& \le & c \|n_+ D^{1/2}_x u_2\|_{L^2_{xt}} \|u_3\|_{L^2_{xt}} \left\langle (\frac{N_2}{N})^\frac{1}{2} 
\right\rangle \\
\nonumber
& \le & c \|n_+\|_{X^{0,\frac{1}{2}+}_+} \|u_2\|_{X^{0+,\frac{1}{2}}} \|u_3\|_{X^{0,0}} \left\langle 
(\frac{N_2}{N})^{\frac{1}{2}} \right\rangle \\
\nonumber
& \le & c \|n_+\|_{X^{0,\frac{1}{2}+}_+} \frac{1}{N_2^{1-}} \|u_2\|_{X^{1,\frac{1}{2}}} \frac{1}{N_3} 
\delta^{\frac{1}{2}-} \|u_3\|_{X^{1,\frac{1}{2}}} \left\langle(\frac{N_2}{N})^{\frac{1}{2}}\right\rangle 
\\
& \le & \frac{c}{N^{1-}} N_{max}^{0-} \delta^{\frac{1}{2}-} \|n_+\|_{X^{0,\frac{1}{2}+}_+} 
\|u_2\|_{X^{1,\frac{1}{2}}} \|u_3\|_{X^{1,\frac{1}{2}}}
\label{9}
\end{eqnarray}
by the bilinear Strichartz estimate. This implies (\ref{8}).

Next we prove (\ref{6}) which is implied by
\begin{eqnarray}
\label{10}
\lefteqn{ \int_0^{\delta} \int_* \left| \frac{m(\xi_2 +\xi_3)-m(\xi_2)m(\xi_3)}{m(\xi_2)m(\xi_3)} \right| 
\widehat{u_1}(\xi_1,t)\widehat{u_2}(\xi_2,t)\widehat{n_+}(\xi_3,t) \, d\xi dt } \\
\nonumber
& \le & \hspace{-0.2cm} c(N^{-\frac{1}{2}+} \delta^{\frac{1}{2}-} + N^{-\frac{3}{2}+}\delta^{0+}) 
\|u_1\|_{X^{-1,\frac{1}{2}}} \|u_2\|_{X^{1,\frac{1}{2}}} \|n_+\|_{X^{0,\frac{1}{2}+}_+}
\end{eqnarray}
{\bf Case 1:} $N_2 \sim N_3 \ge cN $. Then $N_1 \le cN_2$ as above.\\
The multiplier is estimated by $ \frac{c}{m(\xi_2)^2} \le c(\frac{N_2}{N})^{\frac{1}{2}-\epsilon} $, so that 
we get the bound
\begin{eqnarray*}
\lefteqn{ c \|n_+ D^{\frac{1}{2}}_x u_2\|_{L^2_{xt}} \frac{1}{N_2^{\frac{1}{2}}} \|u_1\|_{L^2_{xt}} 
(\frac{N_2}{N})^{\frac{1}{2}-\epsilon} } \\
& \le & c \|n_+\|_{X^{0,\frac{1}{2}+}_+} \|u_2\|_{X^{0+,\frac{1}{2}}} \frac{N_1}{N_2^{\frac{1}{2}}} 
\delta^{\frac{1}{2}-} \|u_1\|_{X^{-1,\frac{1}{2}}} (\frac{N_2}{N})^{\frac{1}{2}-\epsilon} \\
& \le & c \|n_+\|_{X^{0,\frac{1}{2}+}_+} \|u_2\|_{X^{1,\frac{1}{2}}} \frac{N_1}{N_2^{\frac{3}{2}-}} 
\delta^{\frac{1}{2}-} \|u_1\|_{X^{-1,\frac{1}{2}}} (\frac{N_2}{N})^{\frac{1}{2}-\epsilon}
\end{eqnarray*}
which implies (\ref{10}).\\
{\bf Case 2:} $N_1 \sim N_2 \ge cN$, thus $N_3 \le cN_1$.\\
The symbol is majorized by $\frac{c}{m(\xi_2)m(\xi_3)} \le 
c\left\langle(\frac{N_2}{N})^{\frac{1}{2}-\epsilon}\right\rangle$, which can be handled as in Case 1.\\
{\bf Case 3:} $N_1 \sim N_3 \ge cN$ , $N_2 << N_1 \sim N_3$.\\
{\bf Subcase a:} $ N_2 \le N $\\
By the mean value theorem we have
$$ \left| \frac{m(\xi_2 + \xi_3) - m(\xi_2)m(\xi_3)}{m(\xi_2)m(\xi_3)}\right| = \left| \frac{m(\xi_2 + 
\xi_3) - m(\xi_3)}{m(\xi_3)}\right| \le c \left|\frac{(\nabla m)(\xi_3)}{m(\xi_3)}\xi_2\right| \le c 
\frac{N_2}{N_3} $$
and we get the bound
\begin{eqnarray*}
\lefteqn{ c \|n_+ D_x^{\frac{1}{2}}u_1\|_{L^2_{xt}} N_1^{-\frac{1}{2}} \|u_2\|_{L^2_{xt}} \frac{N_2}{N_3} 
} \\
& \le & c \|n_+\|_{X_+^{0,\frac{1}{2}+}} \|u_1\|_{X^{0+,\frac{1}{2}}} N_1^{-\frac{1}{2}} \|u_2\|_{X^{0,0}} 
\frac{N_2}{N_3} \\
& \le & c \|n_+\|_{X_+^{0,\frac{1}{2}+}} \|u_1\|_{X^{-1,\frac{1}{2}}}N_1^{1+} N_1^{-\frac{1}{2}} N_2^{-1} 
\delta^{\frac{1}{2}-} \|u_2\|_{X^{1,\frac{1}{2}}} \frac{N_2}{N_3} \\
& \le & c N_1^{-\frac{1}{2}+} \delta^{\frac{1}{2}-} \|n_+\|_{X_+^{0,\frac{1}{2}+}} 
\|u_1\|_{X^{-1,\frac{1}{2}}} \|u_2\|_{X^{1,\frac{1}{2}}} 
\end{eqnarray*}
which implies (\ref{10}).\\
{\bf Subcase b:} $ |\xi_1| \sim |\xi_3| >> |\xi_2| \ge N $. \\
This is the technically most complicated region where we want to use algebraic manipulations on the 
Fourier side w.r. to $\tau$ and $\xi$ and have also to take into account the characteristic function 
$\psi(t)$ of the time interval $[0,\delta]$. The problem is that $\widehat{\psi}(\tau) = 
\frac{1}{\sqrt{2\pi}} \frac{e^{i\tau \delta}-1}{i \tau} \not\in L^1_{\tau}$, but fortunately $\in 
L_{\tau}^{1+}$. We perform no dyadic decompositions at all here. \\
We estimate the multiplier by $\frac{c}{|m(\xi_2)|} \le c|\xi_2|^{\frac{1}{2}}N^{-\frac{1}{2}}$. Thus our 
aim is to give the following bound
\begin{equation}
 \int_0^\delta \hspace{-0.15cm}\int_* \widehat{u_1}(\xi_1,t) |\xi_2|^{\frac{1}{2}} 
\widehat{u_2}(\xi_2,t)\widehat{n_+}(\xi_3,t) d\xi dt  
 \le c N^{-1+} \delta^{0+} \|u_1\|_{X^{-1,\frac{1}{2}}} \|u_2\|_{X^{1,\frac{1}{2}}} 
\|n_+\|_{X_+^{\frac{1}{2}+}}
\label{***a}
\end{equation}
which would imply (\ref{10}).\\
Abusing notation we denote the Fourier transform w.r. to $x$ and $t$ also by $\; \widehat{} \;$. The 
l.h.s. is bounded by
\begin{equation}
\label{'}
\int_{**} \widehat{u_1}(\xi_1,\tau_1) |\widehat{\psi}(\tau_0)| |\xi_2|^{\frac{1}{2}} 
\widehat{u_2}(\xi_2,\tau_2) \widehat{n_+}(\xi_3,\tau_3) \, d\xi d\tau
\end{equation}
Here ** denotes integration over $ \sum_{i=1}^3 \xi_i = \sum_{i=0}^3 \tau_i = 0 $. Remark again that 
w.l.o.g. $\widehat{u_1},\widehat{u_2},\widehat{n_+} \ge 0 $. The crucial algebraic inequality in our 
region is the following:
$$ |\xi_1| \le c \left( \langle \tau_1 + |\xi_1|^2\rangle^{\frac{1}{2}} + \langle \tau_2 + 
|\xi_2|^2\rangle^{\frac{1}{2}} + \langle \tau_3 + |\xi_3|\rangle^{\frac{1}{2}} + |\tau_0|^{\frac{1}{2}} 
\right) $$
We consider 4 cases according to which of the terms on the r.h.s. is dominant.\\
{\bf Region 1:} $ \langle \tau_1 + |\xi_1|^2 \rangle^{\frac{1}{2}} $ dominant.\\
We get the following bound for (\ref{'}):
\begin{eqnarray*}
\lefteqn{\hspace{-1.7 cm} c \int_{**} \langle \tau_1 + |\xi_1|^2 \rangle^{\frac{1}{2}} |\xi_1|^{-1} 
\widehat{u_1}(\xi_1,\tau_1) |\widehat{\psi}(\tau_0)||\xi_2|^{\frac{1}{2}} \widehat{u_2}(\xi_2,\tau_2) 
\widehat{n_+}(\xi_3,\tau_3)\, d\xi d\tau } \\
& \le & c \|u_1\|_{X^{-1,\frac{1}{2}}}
 \|{\cal F}^{-1}(|\widehat{\psi}|) (D^{1/2}_x u_2) n_+ \|_{L^2_{xt}} \\
& \le & c \|u_1\|_{X^{-1,\frac{1}{2}}} \| {\cal F}^{-1}(|\widehat{\psi}|)\|_{L_t^{\infty -}} \|(D^{1/2}_x 
u_2) n_+\|_{L^{2+}_t L^2_x} \\
& \le & c \delta^{0+} \|u_1\|_{X^{-1,\frac{1}{2}}} \|n_+\|_{X_+^{0,\frac{1}{2}+}} 
\|u_2\|_{X^{0+,\frac{1}{2}}} \\
& \le & c \delta^{0+} N^{-1+} \|u_1\|_{X^{-1,\frac{1}{2}}} \|n_+\|_{X_+^{0,\frac{1}{2}+}} 
\|u_2\|_{X^{1,\frac{1}{2}}}
\end{eqnarray*}
by Lemma \ref{Lemma A'} and by Hausdorff-Young, which gives 
$$ \| {\cal F}^{-1}(|\widehat{\psi}|)\|_{L^{\infty-}_t} \le c \|\widehat{\psi}\|_{L^{1+}_{\tau}} \le c 
\delta^{0+} $$
as one easily calculates.\\
{\bf Region 2:} $ \langle \tau_2 + |\xi_2|^2 \rangle^{\frac{1}{2}} $ dominant.\\
Similarly as before we estimate (\ref{'}) by
\begin{eqnarray*}
\lefteqn{\hspace{-1.7 cm} c \int_{**}  |\xi_1|^{-1} \widehat{u_1}(\xi_1,\tau_1) 
|\widehat{\psi}(\tau_0)|\langle \tau_2 + 
|\xi_2|^2 \rangle^{\frac{1}{2}}|\xi_2|^{\frac{1}{2}} \widehat{u_2}(\xi_2,\tau_2) 
\widehat{n_+}(\xi_3,\tau_3)\, d\xi d\tau } \\
& \le & c \|u_2\|_{X^{\frac{1}{2},\frac{1}{2}}} \|{\cal F}^{-1}(|\widehat{\psi}|) (D^{-1}_x u_1) n_+ 
\|_{L^2_{xt}} \\
& \le & c \|u_2\|_{X^{\frac{1}{2},\frac{1}{2}}} \| {\cal F}^{-1}(|\widehat{\psi}|)\|_{L_t^{\infty -}} 
\|(D^{-1}_x u_1) n_+\|_{L^{2+}_t L^2_x} \\
& \le & c \delta^{0+} \|u_2\|_{X^{\frac{1}{2},\frac{1}{2}}} \|n_+\|_{X_+^{0,\frac{1}{2}+}} 
\|u_1\|_{X^{-\frac{3}{2}+,\frac{1}{2}}} \\
& \le & c \delta^{0+} N^{-\frac{1}{2}} \|u_2\|_{X^{1,\frac{1}{2}}} \|n_+\|_{X_+^{0,\frac{1}{2}+}} 
\|u_1\|_{X^{-1,\frac{1}{2}}} N^{-\frac{1}{2}+}
\end{eqnarray*}
where we used Lemma \ref{Lemma A'} again.\\
{\bf Region 3:} $ \langle \tau_3 + |\xi_3|^2 \rangle^{\frac{1}{2}} $ dominant.\\
Using Lemma \ref{Lemma D}, we control (\ref{'}) by:
\begin{eqnarray*}
\lefteqn{\hspace{-1.7 cm} c \int_{**}  |\xi_1|^{-1} \widehat{u_1}(\xi_1,\tau_1) 
|\widehat{\psi}(\tau_0)| |\xi_2|^{\frac{1}{2}} 
\widehat{u_2}(\xi_2,\tau_2) \langle \tau_3 + |\xi_3|^2 \rangle^{\frac{1}{2}+} 
\widehat{n_+}(\xi_3,\tau_3)\, d\xi d\tau } \\
& \le & c \|n_+\|_{X_+^{0,\frac{1}{2}+}} \|{\cal F}^{-1}(|\widehat{\psi}|) (D^{-1}_x u_1) ( D_x^{\frac{1}{2}} 
u_2) \|_{L^2_{xt}} \\
& \le & c \|n_+\|_{X_+^{0,\frac{1}{2}+}} \| {\cal F}^{-1}(|\widehat{\psi}|)\|_{L_t^{\infty -}} \|(D^{-1}_x 
u_1) ( D_x^{\frac{1}{2}} u_2)\|_{L^{2+}_t L^2_x} \\
& \le & c \delta^{0+} \|n_+\|_{X_+^{0,\frac{1}{2}+}} \|u_1\|_{X^{-\frac{3}{2}+,\frac{1}{2}}} \| 
u_2\|_{X^{\frac{1}{2},\frac{1}{2}}} \\
& \le & c \delta^{0+} \|n_+\|_{X_+^{0,\frac{1}{2}+}} N^{-\frac{1}{2}+} \|u_1\|_{X^{-1,\frac{1}{2}}}  
\|u_2\|_{X^{1,\frac{1}{2}}} N^{-\frac{1}{2}}
\end{eqnarray*}
{\bf Region 4:} $ |\tau_0|^{\frac{1}{2}} $ dominant.\\
The upper bound for (\ref{'}) is here
\begin{eqnarray*}
\lefteqn{\hspace{-1.7 cm} c \int_{**}  |\xi_1|^{-1} \widehat{u_1}(\xi_1,\tau_1) |\tau_0|^{\frac{1}{2}} 
|\widehat{\psi}(\tau_0)| |\xi_2|^{\frac{1}{2}} \widehat{u_2}(\xi_2,\tau_2) \widehat{n_+}(\xi_3,\tau_3)\, 
d\xi d\tau } \\
& \le &  c \|\widehat{D_x^{-1}u_1}\|_{L^2_{\xi_1} L^{1+}_{\tau_1}} \| |\tau|^{\frac{1}{2}} 
|\widehat{\psi}| \ast \widehat{D_x^{1/2} u_2} \ast \widehat{n_+} \|_{L^2_{\xi} L^{\infty -}_{\tau}}
\end{eqnarray*}
by H\"older. The first factor is estimated as follows by H\"older w.r. to $\tau_1$:
\begin{eqnarray*}
\|\widehat{D_x^{-1} u_1}\|_{L^2_{\xi_1} L^{1+}_{\tau_1}} & = & \| \widehat{D_x^{-1} u_1} \, \langle \tau_1 
+ \xi_1^2 \rangle^{\frac{1}{2}} \langle \tau_1 + \xi_1^2 \rangle^{-\frac{1}{2}} \|_{L^2_{\xi_1} 
L^{1+}_{\tau_1}} \\
& \le & \| \widehat{D_x^{-1} u_1} \, \langle \tau_1 + \xi_1^2 \rangle^{\frac{1}{2}}  \|_{L^2_{\xi_1 
\tau_1}} \le c \|u_1\|_{X^{-1,\frac{1}{2}}}
\end{eqnarray*}
The second factor is bounded by Young's inequality by
\begin{eqnarray*}
\lefteqn{ c \| |\tau|^{\frac{1}{2}} |\widehat{\psi}| \|_{L^{2+}_{\tau}} \| \widehat{D^{1/2}_x u_2} \ast 
\widehat{n_+} \|_{L^2_{\xi} L^{2--}_{\tau}} } \\
& \le & c \delta^{0+} \|n_+\|_{X^{0,\frac{1}{2}+}_+} \|u_2\|_{X^{0+,\frac{1}{2}}} \quad \le \quad c 
\delta^{0+} N^{-1+} \|n_+\|_{X^{0,\frac{1}{2}+}_+} \|u_2\|_{X^{1,\frac{1}{2}}}
\end{eqnarray*}
Here we used Lemma \ref{Lemma B} and the bound $ \| |\tau|^{\frac{1}{2}} |\widehat{\psi}| 
\|_{L^{2+}_{\tau}} \le c \delta^{0+} $, which is easily checked.\\
Thus we get (\ref{***a}) in all regions. 

Finally we have to prove (\ref{7}). It is implied by
\begin{eqnarray}
\nonumber
&& \int_0^{\delta} \int_* \left| \frac{m(\xi_1 + \xi_2)}{m(\xi_1)m(\xi_2)} \cdot \frac{m(\xi_3 + \xi_4) - 
m(\xi_3)m(\xi_4)}{m(\xi_3)m(\xi_4)} \right| \widehat{n_+}(\xi_1,t)\widehat{u_2}(\xi_2,t) \cdot \\
&& \nonumber 
\hspace{7cm} \cdot \, \widehat{n_+}(\xi_3,t)\widehat{u_4}(\xi_4,t) d\xi dt  \\
& & \le  c(\frac{1}{N^{3-}} + N^{-1+}\delta^{\frac{1}{2}-}) \|n_+\|_{X^{0,\frac{1}{2}+}_+}^2 
\|u_2\|_{X^{1,\frac{1}{2}}}\|u_4\|_{X^{1,\frac{1}{2}}}
\label{11}
\end{eqnarray}
{\bf Case 1:} $N_1 \sim N_3 \ge cN \, , \, N_1 \sim N_3 >> N_2,N_4 $\\
Subcase a: $N_4 \le N$\\
$$ \left| \frac{m(\xi_1 + \xi_2)}{m(\xi_1)m(\xi_2)} \right| \le \frac{c}{m(\xi_2)} \le c 
\left\langle(\frac{N_2}{N})^{\frac{1}{2}-\epsilon}\right\rangle $$
$$ \left|\frac{m(\xi_3 + \xi_4) - m(\xi_3)m(\xi_4)}{m(\xi_3)m(\xi_4)}\right| = \left|\frac{m(\xi_3 + 
\xi_4) - m(\xi_3)}{m(\xi_3)}\right| \le \left| \frac{(\nabla m)(\xi_3)}{m(\xi_3)} \xi_4 \right| \le 
\frac{cN_4}{N_3} $$
by the mean value theorem. Thus we get the bound
\begin{eqnarray*}
\lefteqn{ c\frac{N_4}{N_3} \|n_+\|_{L^{\infty}_t L^2_x} \|u_2\|_{L^2_t L^{\infty}_x} \|n_+ D^{1/2}_x 
u_4\|_{L^2_{xt}} \frac{1}{N_4^{\frac{1}{2}}} 
\left\langle(\frac{N_2}{N})^{\frac{1}{2}-\epsilon}\right\rangle } \\
& \le & c\frac{N_4}{N_3} \|n_+\|_{X_+^{0,\frac{1}{2}+}} \|u_2\|_{L^2_t H^{\frac{1}{2}+}_x} 
\|n_+\|_{X^{0,\frac{1}{2}+}_+} \| u_4\|_{X^{0+,\frac{1}{2}}} \frac{1}{N_4^{\frac{1}{2}}} 
\left\langle(\frac{N_2}{N})^{\frac{1}{2}-\epsilon}\right\rangle \\
& \le & c\frac{N_4}{N_3} \|n_+\|_{X_+^{0,\frac{1}{2}+}} 
\|u_2\|_{X^{1,\frac{1}{2}}}\frac{1}{N_2^{\frac{1}{2}-}} \delta^{\frac{1}{2}-} 
\|n_+\|_{X^{0,\frac{1}{2}+}_+} \| u_4\|_{X^{1,\frac{1}{2}}} \frac{1}{N_4^{\frac{3}{2}-}} 
\left\langle(\frac{N_2}{N})^{\frac{1}{2}-\epsilon}\right\rangle
\end{eqnarray*}
which implies (\ref{11}).\\
Subcase b: $|\xi_1| \sim |\xi_3| \ge cN \, , \, |\xi_1| \sim |\xi_3| >> |\xi_2|,|\xi_4| \, , \, 
|\xi_2|,|\xi_4| \ge N$.\\
In this case we avoid any dyadic decomposition and estimate as follows:
$$ \left| \frac{m(\xi_1 + \xi_2)}{m(\xi_1)m(\xi_2)} \right| \le \frac{c}{|m(\xi_2)|} \le c 
(\frac{|\xi_2|}{N})^{\frac{1}{2}-} $$
$$ \left|\frac{m(\xi_3 + \xi_4) - m(\xi_3)m(\xi_4)}{m(\xi_3)m(\xi_4)}\right| \le \frac{c}{|m(\xi_4)|} \le 
c(\frac{|\xi_4|}{N})^{\frac{1}{2}-}  $$
Thus we get the bound
\begin{eqnarray*}
&& \hspace{-1cm}\int_0^{\delta} \int_* \frac{|\xi_2|^{\frac{1}{2}-}}{N^{\frac{1}{2}-}} 
\frac{|\xi_4|^{\frac{1}{2}-}}{N^{\frac{1}{2}-}} 
\widehat{n_+}(\xi_1,t)(|\xi_2|^{\frac{1}{2}-}\widehat{u_2}(\xi_2,t))\frac{1}{|\xi_2|^{\frac{1}{2}-}} \cdot 
\\
&& \cdot \, 
\widehat{n_+}(\xi_3,t)(|\xi_4|^{\frac{1}{2}-}\widehat{u_4}(\xi_4,t))\frac{1}{|\xi_4|^{\frac{1}{2}-}} \, 
d\xi dt \\
& \le & \frac{c}{N^{1-}} \|n_+\|_{L^{\infty}_t L^2_x} \|D^{\frac{1}{2}-}_x u_2\|_{L^2_t L^{\infty}_x} 
\|n_+\|_{L^{\infty}_t L^2_x} \|D^{\frac{1}{2}-}_x u_4\|_{L^2_t L^{\infty}_x} \\
& \le & \frac{c}{N^{1-}} \|n_+\|_{X_+^{0,\frac{1}{2}+}} \|u_2\|_{X^{1,\frac{1}{2}}} \delta^{\frac{1}{2}-} 
\|n_+\|_{X_+^{0,\frac{1}{2}+}} \|u_4\|_{X^{1,\frac{1}{2}}} \delta^{\frac{1}{2}-}
\end{eqnarray*}
which is sufficient, because no dyadic decomposition was performed.\\
Subcase c: $|\xi_1| \sim |\xi_3| \ge cN \, , \, |\xi_1| \sim |\xi_3| >> |\xi_2|,|\xi_4| \, , \, |\xi_4| 
\ge N \ge |\xi_2|$.\\
We again perform no dyadic decomposition and estimate as follows:
$$ \left| \frac{m(\xi_1 + \xi_2)}{m(\xi_1)m(\xi_2)} \right| \le c  $$
$$ \left|\frac{m(\xi_3 + \xi_4) - m(\xi_3)m(\xi_4)}{m(\xi_3)m(\xi_4)}\right| \le \frac{c}{|m(\xi_4)|} \le 
c(\frac{|\xi_4|}{N})^{\frac{1}{2}}  $$
Thus we get the bound
\begin{eqnarray*}
\lefteqn{ \int_0^{\delta} \int_*  \frac{|\xi_4|^{\frac{1}{2}}}{N^{\frac{1}{2}}} 
\widehat{n_+}(\xi_1,t)\widehat{u_2}(\xi_2,t)) 
\widehat{n_+}(\xi_3,t)(|\xi_4|^{\frac{1}{2}}\widehat{u_4}(\xi_4,t))\frac{1}{|\xi_4|^{\frac{1}{2}}} \, d\xi 
dt }\\
& \le & \frac{c}{N^{\frac{1}{2}}} \|n_+\|_{L^{\infty}_t L^2_x} \|u_2\|_{L^2_t L^{\infty}_x} \|n_+ 
D^{1/2}_x u_4\|_{L^2_{xt}} \\
& \le & \frac{c}{N^{\frac{1}{2}}} \|n_+\|_{X^{0,\frac{1}{2}+}_+} \|u_2\|_{L^2_t H^{\frac{1}{2}+}_x} 
\|n_+\|_{X^{0,\frac{1}{2}+}_+} \| u_4\|_{X^{0+,\frac{1}{2}}} \\
& \le & \frac{c}{N^{\frac{1}{2}}} \|n_+\|_{X^{0,\frac{1}{2}+}_+} 
\delta^{\frac{1}{2}-}\|u_2\|_{X^{1,\frac{1}{2}}} \|n_+\|_{X^{0,\frac{1}{2}+}_+} \frac{1}{N^{1-}} \| 
u_4\|_{X^{1,\frac{1}{2}}} \\
& \le & \frac{c}{N^{\frac{3}{2}-}} \delta^{\frac{1}{2}-} \|n_+\|_{X^{0,\frac{1}{2}+}_+}^2 
\|u_2\|_{X^{1,\frac{1}{2}}} \| u_4\|_{X^{1,\frac{1}{2}}}
\end{eqnarray*}
{\bf Case 2:} $N_2 \sim N_4 \ge cN \, , \, N_2 \sim N_4 >> N_1,N_3 $ \\
We have
$$ \left| \frac{m(\xi_1 + \xi_2)}{m(\xi_1)m(\xi_2)} \right| \le \frac{c}{|m(\xi_1)|}  \le c 
\left\langle(\frac{N_1}{N})^{\frac{1}{2}}\right\rangle \le c (\frac{N_4}{N})^{\frac{1}{2}}$$
$$ \left|\frac{m(\xi_3 + \xi_4) - m(\xi_3)m(\xi_4)}{m(\xi_3)m(\xi_4)}\right|  \le \frac{c}{|m(\xi_3)|}  
\le c \left\langle(\frac{N_3}{N})^{\frac{1}{2}}\right\rangle \le c (\frac{N_2}{N})^{\frac{1}{2}}$$
This gives the bound
\begin{eqnarray*}
\lefteqn{ c \|n_+ D^{1/2}_x u_2\|_{L^2_{xt}} N_2^{-1/2} \|n_+ D^{1/2}_x u_4\|_{L^2_{xt}} N_4^{-1/2} 
\left(\frac{N_4}{N}\right)^{1/2} \left(\frac{N_2}{N}\right)^{1/2} } \\
& \le & \frac{c}{N} \|n_+\|_{X^{0,\frac{1}{2}+}_+} \|u_2\|_{X^{0+,\frac{1}{2}}} 
\|n_+\|_{X^{0,\frac{1}{2}+}_+} \|u_4\|_{X^{0+,\frac{1}{2}}} \\
& \le & \frac{c}{N} N_2^{-1+} N_4^{-1+} \|n_+\|_{X^{0,\frac{1}{2}+}_+}^2 \|u_2\|_{X^{1,\frac{1}{2}}}  
\|u_4\|_{X^{1,\frac{1}{2}}} 
\end{eqnarray*}
which implies (\ref{11}).\\
{\bf Case 3:} $N_1 \sim N_2 \ge cN ß, , \, N_1 \sim N_2 >> N_3,N_4 $\\
Using
$$ \left| \frac{m(\xi_1 + \xi_2)}{m(\xi_1)m(\xi_2)} \right| \le \frac{c}{|m(\xi_1)|^2}  \le c 
(\frac{N_1}{N})^{\frac{1}{2}} $$
$$ \left|\frac{m(\xi_3 + \xi_4) - m(\xi_3)m(\xi_4)}{m(\xi_3)m(\xi_4)}\right|  \le 
\frac{c}{|m(\xi_3)m(\xi_4)|}  \le c \left\langle(\frac{N_3}{N})^{\frac{1}{2}}\right\rangle 
\left\langle(\frac{N_4}{N})^{\frac{1}{2}}\right\rangle$$
we get the bound
\begin{eqnarray*}
\lefteqn{ c \|n_+ D^{1/2}_x u_2\|_{L^2_{xt}} N_2^{-1/2} \|n_+\|_{L^{\infty}_t L^2_x} \|u_4\|_{L^2_t 
L^{\infty}_x} \left\langle (\frac{N_3}{N})^{\frac{1}{2}}\right\rangle \left\langle 
(\frac{N_4}{N})^{\frac{1}{2}}\right\rangle (\frac{N_1}{N})^{\frac{1}{2}} } \\
 & \hspace{-0.2cm} \le & \hspace{-0.2cm} c \|n_+\|_{X^{0,\frac{1}{2}+}_+} \|u_2\|_{X^{0+,\frac{1}{2}}} 
N_2^{-1/2} \|n_+\|_{X^{0,\frac{1}{2}+}_+} \|u_4\|_{L^2_t H^{\frac{1}{2}+}_x} \left\langle 
(\frac{N_3}{N})^{\frac{1}{2}}\right\rangle \left\langle (\frac{N_4}{N})^{\frac{1}{2}}\right\rangle 
(\frac{N_1}{N})^{\frac{1}{2}} 
\\
& \hspace{-0.2cm} \le & \hspace{-0.2cm} c \|n_+\|_{X^{0,\frac{1}{2}+}_+} \|u_2\|_{X^{1,\frac{1}{2}}} 
N_2^{-\frac{3}{2}+} \|n_+\|_{X^{0,\frac{1}{2}+}_+} \|u_4\|_{X^{1,\frac{1}{2}}}\delta^{\frac{1}{2}-} 
N_4^{-\frac{1}{2}+} \cdot \\ && \hspace{7cm}\cdot \,  \left\langle 
(\frac{N_2}{N})^{\frac{1}{2}}\right\rangle \left\langle (\frac{N_4}{N})^{\frac{1}{2}}\right\rangle 
(\frac{N_1}{N})^{\frac{1}{2}} \\
& \hspace{-0.2cm} \le & \hspace{-0.2cm} c N^{-\frac{3}{2}+} (N_1N_2N_3N_4)^{0-} \delta^{\frac{1}{2}-} 
\|n_+\|_{X^{0,\frac{1}{2}+}_+}^2 \|u_2\|_{X^{1,\frac{1}{2}}} \|u_4\|_{X^{1,\frac{1}{2}}}
\end{eqnarray*}
{\bf Case 4:} $N_2 \sim N_3 \ge cN \, , \, N_2 \sim N_3 >> N_1,N_4$ \\
Using
$$ \left| \frac{m(\xi_1 + \xi_2)}{m(\xi_1)m(\xi_2)} \right| \le \frac{c}{|m(\xi_1)|}  \le c 
\left\langle(\frac{N_1}{N})^{\frac{1}{2}}\right\rangle $$
$$ \left|\frac{m(\xi_3 + \xi_4) - m(\xi_3)m(\xi_4)}{m(\xi_3)m(\xi_4)}\right|  \le \frac{c}{|m(\xi_4)|}  
\le c  \left\langle(\frac{N_4}{N})^{\frac{1}{2}}\right\rangle$$
we get the bound
\begin{eqnarray*}
\lefteqn{ c \|n_+ D^{1/2}_x u_2\|_{L^2_{xt}} N_2^{-1/2} \|n_+\|_{L^{\infty}_t L^2_x} \|u_4\|_{L^2_t 
L^{\infty}_x} \left\langle (\frac{N_1}{N})^{\frac{1}{2}}\right\rangle \left\langle 
(\frac{N_4}{N})^{\frac{1}{2}}\right\rangle } \\
&  \le &  c \|n_+\|_{X^{0,\frac{1}{2}+}_+} \|u_2\|_{X^{1,\frac{1}{2}}} N_2^{-\frac{3}{2}+} 
\|n_+\|_{X^{0,\frac{1}{2}+}_+} \|u_4\|_{X^{1,\frac{1}{2}}}\delta^{\frac{1}{2}-} N_4^{-\frac{1}{2}+}
 \cdot \\ && \hspace{7cm}
 \cdot \,  \left\langle (\frac{N_2}{N})^{\frac{1}{2}}\right\rangle \left\langle 
(\frac{N_4}{N})^{\frac{1}{2}}\right\rangle \\
&  \le &  c N^{-\frac{3}{2}+} (N_1N_2N_3N_4)^{0-} \delta^{\frac{1}{2}-} \|n_+\|_{X^{0,\frac{1}{2}+}_+}^2 
\|u_2\|_{X^{1,\frac{1}{2}}} \|u_4\|_{X^{1,\frac{1}{2}}}
\end{eqnarray*}
{\bf Case 5:} $N_3 \sim N_4\ge cN \, , \, N_3\sim N_4>> N_1,N_2$\\
Using
$$ \left| \frac{m(\xi_1 + \xi_2)}{m(\xi_1)m(\xi_2)} \right| \le \frac{c}{|m(\xi_1)m(\xi_2)|}  \le c 
\left\langle(\frac{N_1}{N})^{\frac{1}{2}}\right\rangle 
\left\langle(\frac{N_2}{N})^{\frac{1}{2}}\right\rangle $$
$$ \left|\frac{m(\xi_3 + \xi_4) - m(\xi_3)m(\xi_4)}{m(\xi_3)m(\xi_4)}\right|  \le \frac{c}{|m(\xi_4)|^2}  
\le c (\frac{N_4}{N})^{\frac{1}{2}}$$
we get the bound
\begin{eqnarray*}
\lefteqn{ c  \|n_+\|_{L^{\infty}_t L^2_x} \|u_2\|_{L^2_t L^{\infty}_x} \|n_+ D^{1/2}_x u_4\|_{L^2_{xt}} 
N_4^{-1/2}\left\langle (\frac{N_1}{N})^{\frac{1}{2}}\right\rangle \left\langle 
(\frac{N_2}{N})^{\frac{1}{2}}\right\rangle (\frac{N_4}{N})^{\frac{1}{2}} } \\
&  \le &  c \|n_+\|_{X^{0,\frac{1}{2}+}_+} \|u_2\|_{X^{1,\frac{1}{2}}} 
N_2^{-\frac{1}{2}+}\delta^{\frac{1}{2}-} \|n_+\|_{X^{0,\frac{1}{2}+}_+} \|u_4\|_{X^{1,\frac{1}{2}}} 
N_4^{-\frac{3}{2}+}
 \cdot \\ && \hspace{6cm}
 \cdot \, \left\langle (\frac{N_1}{N})^{\frac{1}{2}}\right\rangle \left\langle 
(\frac{N_2}{N})^{\frac{1}{2}}\right\rangle (\frac{N_4}{N})^{\frac{1}{2}} \\
&  \le &  c N^{-\frac{3}{2}+} (N_1N_2N_3N_4)^{0-} \delta^{\frac{1}{2}-} \|n_+\|_{X^{0,\frac{1}{2}+}_+}^2 
\|u_2\|_{X^{1,\frac{1}{2}}} \|u_4\|_{X^{1,\frac{1}{2}}}
\end{eqnarray*}
{\bf Case 6:} $N_1 \sim N_4\ge cN \, , \, N_1\sim N_4>> N_2,N_3$\\
Using
$$ \left| \frac{m(\xi_1 + \xi_2)}{m(\xi_1)m(\xi_2)} \right| \le \frac{c}{|m(\xi_2)|}  \le c 
\left\langle(\frac{N_2}{N})^{\frac{1}{2}}\right\rangle $$
$$ \left|\frac{m(\xi_3 + \xi_4) - m(\xi_3)m(\xi_4)}{m(\xi_3)m(\xi_4)}\right|  \le \frac{c}{|m(\xi_3)|}  
\le c \left\langle(\frac{N_3}{N})^{\frac{1}{2}}\right\rangle $$
we get the bound
\begin{eqnarray*}
\lefteqn{ c  \|n_+\|_{L^{\infty}_t L^2_x} \|u_2\|_{L^2_t L^{\infty}_x} \|n_+ D^{1/2}_x u_4\|_{L^2_{xt}} 
N_4^{-1/2}\left\langle (\frac{N_2}{N})^{\frac{1}{2}}\right\rangle \left\langle 
(\frac{N_3}{N})^{\frac{1}{2}}\right\rangle } \\
&  \le &  c \|n_+\|_{X^{0,\frac{1}{2}+}_+} \|u_2\|_{X^{1,\frac{1}{2}}} 
N_2^{-\frac{1}{2}+}\delta^{\frac{1}{2}-} \|n_+\|_{X^{0,\frac{1}{2}+}_+} \|u_4\|_{X^{1,\frac{1}{2}}} 
N_4^{-\frac{3}{2}+}
 \cdot \\ && \hspace{6cm}
 \cdot \, \left\langle (\frac{N_2}{N})^{\frac{1}{2}}\right\rangle \left\langle 
(\frac{N_3}{N})^{\frac{1}{2}}\right\rangle \\
&  \le &  c N^{-\frac{3}{2}+} (N_1N_2N_3N_4)^{0-} \delta^{\frac{1}{2}-} \|n_+\|_{X^{0,\frac{1}{2}+}_+}^2 
\|u_2\|_{X^{1,\frac{1}{2}}} \|u_4\|_{X^{1,\frac{1}{2}}}
\end{eqnarray*}
The remaining cases where at least three factors have equivalent frequencies $\ge cN$ are similar or 
easier to handle so that (\ref{11}) is proved in all possible situations. \\
The proof of the proposition is complete.

\section{The global existence result}
\begin{theorem}
Let $ 1 > s > 5/6 $. The Zakharov system (\ref{1}),(\ref{2}),(\ref{3}) has a unique global solution for 
data $ u_0 \in H^s({\bf R})$ , $n_0 \in H^{s-1}({\bf R}) $ , $ A^{-1/2}n_1 \in  H^{s-1}({\bf R}) $. More 
precisely, for any $T > 0$ there exists a unique solution
$$ (u,n,A^{-1/2}n_t) \in X^{s,\frac{1}{2}}[0,T] \times \widetilde{X}^{s-1,\frac{1}{2}+}[0,T] \times  
\widetilde{X}^{s-1,\frac{1}{2}+}[0,T] $$
where $ \widetilde{X}^{s-1,\frac{1}{2}+}[0,T] := X_+^{s-1,\frac{1}{2}+}[0,T] + X_-^{s-1,\frac{1}{2}+}[0,T] $. This solution satisfies
$$ (u,n,A^{-1/2}n_t) \in C^0([0,T],H^s({\bf R}) \times H^{s-1}({\bf R}) \times  H^{s-1}({\bf R})) $$
and
$$ \|u(t)\|_{H^s} + \|n(t)\|_{H^{s-1}} + \|A^{-1/2} n_t(t)\|_{H^{s-1}} \le c (1+t)^{\frac{2(1-s)}{6s-5}+} 
$$
\end{theorem}
{\bf Proof:} The data satisfy the estimates
\begin{eqnarray*}
\|Iu_0\|_{H^1} & \le & c N^{1-s}\|u_0\|_{H^s} \\
\|In_{\pm 0}\|_{L^2} & \le & c N^{1-s} (\|n_0\|_{H^{s-1}} +\|A^{-1/2} n_1\|_{H^{s-1}})
\end{eqnarray*}
We use our local existence theorem on $[0,\delta]$, where $ \delta \sim \frac{1}{N^{4(1-s)+}} $ and 
conclude
\begin{eqnarray}
\nonumber
\lefteqn{ \|Iu\|_{X^{1,\frac{1}{2}}[0,\delta]} + \|In_+\|_{X^{0,\frac{1}{2}+}_+[0,\delta]} + 
\|In_-\|_{X^{0,\frac{1}{2}+}_-[0,\delta]} } \\
&& \le  c(\|Iu_0\|_{H^1} + \|In_+\|_{L^2} + \|In_-\|_{L^2}) \le c_2 N^{1-s}
\label{**}
\end{eqnarray}
From (\ref{*'}) we get
$$ E(Iu_0,In_{+0}) \le c_0 (\|Iu_0\|_{H^1}^2 + \|In_{+0}\|_{L^2}^2 + \|Iu_0\|_{L^2}^6) \le \overline{c} 
N^{2(1-s)} $$
and from (\ref{*})
$$ \|A^{1/2}Iu_0\|_{L^2}^2 + \|In_+\|_{L^2}^2 + \|In_-\|_{L^2}^2 \le \widehat{c} N^{2(1-s)} \quad , \quad 
\|Iu_0\|_{L^2} \le M $$
with $ \widehat{c} = \widehat{c}(\overline{c}) $. Thus the constant in (\ref{**}) depends only on 
$\overline{c}$ and $M$, i.e. $c_2 = c_2(\overline{c},M)$.\\
In order to reapply the local existence result with time intervals of equal length we need a uniform bound 
of the solution at time $t=\delta$ and $t=2\delta$ etc. which follows from a uniform control over the 
energy by (\ref{*}). The increment of the energy is controlled by Proposition \ref{Prop. 3.1} and 
(\ref{**}) as follows:
\begin{eqnarray*}
\lefteqn{ |E(Iu(\delta),In_+(\delta)) - E(Iu_0,In_{+0})| } \\
& \le & c[(N^{-\frac{1}{2}+} \delta^{\frac{1}{2}-} + N^{-\frac{3}{2}+}) 
\|In_+\|_{X^{0,\frac{1}{2}+}_+[0,\delta]} \|Iu\|_{X^{1,\frac{1}{2}}[0,\delta]}^2 \\
&& + (N^{-3+} + N^{-1+} \delta^{\frac{1}{2}-})  \|In_+\|_{X^{0,\frac{1}{2}+}_+[0,\delta]}^2 
\|Iu\|_{X^{1,\frac{1}{2}}[0,\delta]}^2] \\
& \le & c((N^{-\frac{1}{2}+} \delta^{\frac{1}{2}-} + N^{-\frac{3}{2}+}) N^{3(1-s)} +(N^{-3+} + N^{-1+} 
\delta^{\frac{1}{2}-}) N^{4(1-s)})
\end{eqnarray*}
Using the definition of $\delta$ we arrive at
\begin{eqnarray*}
\lefteqn{ |E(Iu(\delta),In_+(\delta)) - E(Iu_0,In_{+0})| } \\
& \hspace{-0.2 cm} \le & \hspace{-0.2 cm} c_3 ((N^{-\frac{1}{2}+} N^{-2(1-s)+}+ N^{-\frac{3}{2}+}) 
N^{3(1-s)} + (N^{-3+} + N^{-1+} N^{-2(1-s)+})N^{4(1-s)}) \\
& \hspace{-0.2 cm} \le & \hspace{-0.2 cm} c_3( N^{-\frac{1}{2}+} N^{-2(1-s)+} N^{3(1-s)} + N^{-1+} 
N^{-2(1-s)+} N^{4(1-s)})
\end{eqnarray*}
where $ c_3 = c_3(\overline{c},M) $. This is easily seen to be bounded by $ \overline{c}N^{2(1-s)}$ (for 
large $N$).\\
The number of iteration steps to reach the given time $T$ is $\frac{T}{\delta} \sim T N^{4(1-s)+} $. This 
means that in order to give a uniform bound of the energy of the iterated solutions, namely by 
$2\overline{c}N^{2(1-s)}$, from the last inequality the following condition has to be fulfilled:
$$
c_3(N^{-\frac{1}{2}+}N^{-2(1-s)}N^{3(1-s)} + N^{-1+}N^{-2(1-s)+}N^{4(1-s)})TN^{4(1-s)+}   < \overline{c} 
N^{2(1-s)} $$ 
where $c_3 = c_3(2\overline{c},2M)$ (recall here that the initial energy is bounded by 
$\overline{c}N^{2(1-s)}$).\\
This can be fulfilled for $N$ sufficiently large provided the following conditions hold:
\begin{eqnarray*}
-\frac{1}{2}-2(1-s)+3(1-s)+4(1-s)<2(1-s) & \Longleftrightarrow & s > 5/6 \\
-1-2(1-s)+4(1-s)+4(1-s)<2(1-s) & \Longleftrightarrow & s > 3/4
\end{eqnarray*}
So here is the point where the decisive bound on $s$ appears.\\
A uniform bound of the energy implies by (\ref{*}) uniform control of
$$ \|A^{1/2}Iu(t)\| + \|In(t)\| + \|A^{-1/2}In_t(t)\| \le c N^{1-s}$$
Moreover $ \|Iu(t)\| \le \|u(t)\| = \|u_0\| $ , thus
$$ \|u(t)\|_{H^s} + \|n(t)\|_{H^{s-1}} + \|A^{-1/2}n_t(t)\|_{H^{s-1}} \le c N^{1-s} $$
Now, one can directly give a bound on the growth of the solution as follows. The most restrictive 
condition on $N$ comes from the inequality
$$c_3 T N^{-\frac{1}{2}+} N^{-2(1-s)} N^{3(1-s)} N^{4(1-s)+} < \overline{c} N^{2(1-s)} \Longleftrightarrow 
N > c T^{\frac{2}{6s-5}+} $$
This implies
$$ \sup_{0\le t \le T} (\|u(t)\|_{H^s} + \|n(t)\|_{H^{s-1}} + \|A^{-1/2}n_t(t)\|_{H^{s-1}}) \le c 
(1+T)^{\frac{2(1-s)}{6s-5}+} $$


\begin{thebibliography}{AAAAA1}
\bibitem[B1]{B1} J.  Bourgain: {\sl Refinements of Strichartz' inequality and applications to 2D-NLS with 
critical nonlinearity}.
Int. Math. Res. Not. no. 5 (1998), 253-283
\bibitem[B2]{B2} J.  Bourgain: {\sl Scattering in the energy space and below for 3D NLS}.  
J. d'Anal. Math. 75 (1998), 267-297 
\bibitem[B3]{B3} J.  Bourgain: {\sl Global solutions of nonlinear Schr\"odinger equations}.  
Amer. Math. Soc. Colloq. Publ. 46, Amer. Math. Soc., Providence, 1999 
\bibitem[CKSTT1]{CKSTT1} J. Colliander, M. Keel, G. Staffilani, H. Takaoka, and T. Tao: {\sl Global 
well-posedness for Schr\"odinger equations with derivative}.
Siam J. Math. Analysis 33 (2001), 649-669
\bibitem[CKSTT2]{CKSTT2} J. Colliander, M. Keel, G. Staffilani, H. Takaoka, and T. Tao: {\sl
A refined global well-posedness result for Schr\"odinger equations with derivative}. Siam J. Math. 
Analysis 34 (2002), 64-86
\bibitem[CKSTT3]{CKSTT3} J. Colliander, M. Keel, G. Staffilani, H. Takaoka, and T. Tao: {\sl
Almost conservation laws and global rough solutions to a nonlinear Schr\"odinger
equation}. Math. Res. Letters 9 (2002), 659-682
\bibitem[CKSTT4]{CKSTT4} J. Colliander, M. Keel, G. Staffilani, H. Takaoka, and T. Tao: {\sl
Sharp global well-posedness for KdV and modified KdV on {\bf R} and {\bf T}}. J. Amer. Math. Soc. 16 
(2003), 705-749
\bibitem[CKSTT5]{CKSTT5} J. Colliander, M. Keel, G. Staffilani, H. Takaoka, and T. Tao: {\sl Multilinear 
estimates for periodic KdV equations, and applications}.
J. Funct. Anal. 211 (2004), 173-218
\bibitem[CKSTT6]{CKSTT6} J. Colliander, M. Keel, G. Staffilani, H. Takaoka, and T. Tao: {\sl
Global existence and scattering for rough solutions to a nonlinear Schr\"odinger equation on ${\bf R}^3$}. 
Comm. Pure Appl. Math. 57 (2004), 987-1014
\bibitem[CST]{CST} J. Colliander, G. Staffilani, and H. Takaoka: {\sl Global wellposedness for KdV below 
$L^2$}. Math. Res. Lett. 6 (1999), 755-778 
\bibitem[FLP]{FLP} G. Fonseca, F. Linares, and G. Ponce: {\sl Global well-posedness for the modified 
Korteweg-de Vries equation}. Comm. Partial Differential Equations 24 (1999), 683-705
\bibitem[GTV]{GTV} J. Ginibre, Y.  Tsutsumi and G.  Velo: {\sl On the Cauchy problem for the Zakharov 
system}. J. Funct. Anal. 151 (1997), 384-436 
\bibitem[G]{G} A. Gr\"unrock: {\sl New applications of the Fourier restriction norm method to 
wellposedness problems for nonlinear evolution equations}. Dissertation Univ. Wuppertal 2002
\bibitem[KT]{KT} M.  Keel, and T.  Tao: {\sl Local and global well-posedness of wave maps on ${\bf 
R}^{1+1}$ for rough data}. Int. Math. Res. Not. no. 21 (1998), 1117-1156
\bibitem[KPV]{KPV} C.  Kenig, G.  Ponce and L.  Vega:  {\sl Global well-posedness for semi-linear wave 
equations}. Comm. Part. Diff. Equ. 25 (2000), 1741-1752 
\bibitem[KS]{KS} S. Klainerman and S. Selberg: {\sl Bilinear estimates and applications to nonlinear wave 
equations}. Comm. Contemp. Math. 4 (2002), 223-295
\bibitem[P1]{P1} H. Pecher: {\sl Global well-posedness below energy space for the 1-dimensional Zakharov 
system}. Int. Math. Res. Not. no. 19 (2001), 1027-1056
\bibitem[P2]{P2} H. Pecher: {\sl Global solutions of the Klein-Gordon-Schr\"odinger 
system with rough data}. Differential Integral Equations 17 (2004), 179-214
\bibitem[Z]{Z} Zakharov, V.E.: {\sl Collapse of Langmuir waves}. Soviet Phys. JETP 35 (1972), 908-914 
\end{thebibliography}
\end{document}